\documentclass[12pt]{article}
\usepackage{amsmath}
\usepackage{amssymb}
\usepackage{latexsym}
\usepackage{graphicx}

\title
{A spectral analogue of the Meinardus theorem on asymptotics of 
the number of partitions}
\author
{
Tatsuya Tate\thanks{Research partially supported by 
JSPS Grant-in-Aid for Scientific Research (No.~18740089, No.~21740117).}\\
Graduate School of Mathematics \\
Nagoya University \\
Furo-cho, Chikusa-ku, \\
Nagoya, 464--8602, Japan\\
Email: tate@math.nagoya-u.ac.jp}
\date{\empty}

\newcommand{\vol}{{\operatorname{Vol}}}

\renewcommand{\phi}{\varphi}

\newcommand{\ispa}[1]{\langle \,#1 \,\rangle }

\newcommand{\sgn}{\operatorname{sgn}\nolimits} 
\newcommand{\tr}{\operatorname{Tr}\nolimits}

\newcommand{\mb}{\mathbb}

\newcommand{\hcal}{\mathcal{H}}

\newcommand{\ocal}{\mathcal{O}}

\newcommand{\ep}{\varepsilon}

\newcommand{\re}{{\rm Re}\,}
\newcommand{\im}{{\rm Im}\,}

\newtheorem{theo}{{\sc Theorem}}[section]
\newtheorem{maintheo}{{\sc Theorem}}

\newtheorem{maincor}[maintheo]{{\sc Corollary}}
\newtheorem{lem}[theo]{{\sc Lemma}}
\newtheorem{prop}[theo]{{\sc Proposition}}

\newenvironment{proof}%
{\def\psymbol{{\it Proof. }\enspace}
\psymbol}%
{\def\qed{\rule[-2pt]{3pt}{9pt}}
\hspace{3pt}\qed\par\bigskip}

\begin{document}
\maketitle

\renewcommand{\labelenumi}{{\rm (\arabic{enumi})}}
\maketitle

\begin{abstract}
We discuss asymptotics of the number of states of Boson gas 
whose Hamiltonian is given by a positive elliptic pseudo-differential operator 
of order one on a compact manifold. 
We obtain an asymptotic formula for the average of the number of states. 
Furthermore, 
when the operator has integer eigenvalues and the periodic 
orbits of period less than $2\pi$ of the classical mechanics form 
clean submanifolds of lower dimensions, we give an asymptotic formula 
for the number of states itself.  
This is regarded as an analogue of the Meinardus theorem 
on asymptotics of the number of partitions of a positive integer. 
We use the Meinardus saddle point method of obtaining the 
asymptotics of the number of partitions, combined with a theorem due to 
Duistermaat-Guillemin and other authors on the singularities 
of the trace of the wave operators. 
\end{abstract}

\section{Introduction}
\label{Intro}

The subject on asymptotics of the number of partitions $p(E)$ of a positive integer $E$ 
is one of the main topics in classical analytic number theory, and there are many results in this area 
since the celebrated work of Hardy-Ramanujan (\cite{HR2}) appeared. As a simplest 
form of their theorem, one has 
\begin{equation}
\label{hra}
p(E) \sim \frac{1}{4E\sqrt{3}}e^{\pi 
\left(
\frac{2E}{3}
\right)^{1/2}},\qquad E \to \infty. 
\end{equation}
There are many other works on the asymptotics of $p(E)$ or its variants. 
For example, after Hardy-Ramanujan's work, Rademacher (\cite{Ra}) gives an 
expansion of $p(E)$ in a convergent series. See \cite{An}, \cite{EL}, \cite{Ri1}, \cite{Ri2} 
for some of them. 

The problem which we are going to address here is on asymptotic 
behavior of the number of states of Boson gas which is, as in the title, 
regarded as a spectral analogue of the number of partitions. 

In this section, we explain the problem which we consider in this paper, 
and we state the main results after  we give some accounts on the motivation for the problem. 
The method we use to obtain asymptotics of the number of states of Boson gas 
is Meinardus' suddle point method (\cite{An}, \cite{Me}) for obtaining asymptotic 
formula of the number of partitions. This method has some applications to mathematical physics, 
such as finding asymptotic growth of the state density for $p$-branes (\cite{Ac}, \cite{BKZ}, \cite{EORBZ}). 
We review, also in this section, Meinardus' theorem and explain roughly the method.

\subsection{A problem in quantum statistics}
\label{Intro1}

Let $P$ be a first order strictly positive elliptic pseudo-differential operator of order one 
on a compact connected smooth manifold $M$ of dimension $n$ with 
principal symbol $0<p \in C^{\infty}(T^{*}M \setminus 0)$. 
For each positive integer $N$, consider the operator 
\begin{equation}
\label{gop}
P^{N}=\sum_{i=1}^{N}
I \otimes \cdots \otimes  I \otimes 
\operatornamewithlimits{{\it P}}^{(i)} \otimes I 
\otimes \cdots \otimes I
\end{equation}
on the $L^{2}$-space $L^{2}(M^{N})$ on the $N$-fold product space 
$M^{N}=M \times \cdots \times M$ ($N$-times). 
We have the natural action of the symmetric 
goup $\mathfrak{S}_{N}$ on $L^{2}(M^{N})$, and 
this action commutes with the operator $P^{N}$. 
Therefore, the operator $P^{N}$ can be restricted to the symmetric and anti-symmetric subspaces 
$\hcal_{S}^{N}$, $\hcal_{A}^{N}$ (see Section \ref{Partition} for the precise definition). 
The restrictions  $P^{N}_{B}=P^{N}|_{\hcal_{S}^{N}}$ and $P^{N}_{F}=P^{N}|_{\hcal_{A}^{N}}$ 
are regarded as the Hamiltonian for the system of $N$-particles consisting of Boson and Fermion, respectively, 
each particle of which is dominated by the Hamiltonian $P$, and it is assumed that there are no 
interaction among each particle. 
In this paper, we restrict our attention to the systems of Boson. 
Let ${\rm Spec}(P^{N}_{B})$ denote the spectrum of the Bose Hamiltonian, $P^{N}_{B}$. 
It is easy to see that ${\rm Spec}(P^{N}_{B})={\rm Spec}(P^{N})$, and hence ${\rm Spec}(P^{N}_{B})$ consists 
of eigenvalues of finite multiplicities. 
But the multiplicities of each eigenvalues of the Bose Hamiltonian $P^{N}_{B}$ are different from that of 
$P^{N}$. Denote the multiplicity of the eigenvalue $E \in {\rm Spec}(P^{N}_{B})$ by $\Omega(N,E)$. 
The quantity $\Omega (N,E)$ is naturally regarded as a number of states of $N$-Boson system 
in the energy $E$.

Natural question here is to find asymptotic formulas for $\Omega(N,E)$ when 
$E$ and $N$ get large and $E$ and $N$ are related by a relation, like $E/N \sim {\rm const}$. 
But, it seems to be difficult to find them. 
Thus, it would be also natural to consider the asymptotics as $N \to \infty$ 
of the density of states, which in this case would be given by 
\[
\sum_{E' \in {\rm Spec}(P^{N}_{B})}\Omega (N,E')\delta_{E'}
\ \ \ \mbox{or more specifically} \ \ \sum_{E' \in [E,E+e_{N}]}\Omega (N,E')
\]
with a suitable sequence $e_{N}>0$, however it is not so transparent how to choose $e_{N}$. 
These two problems would be rather natural in physical points of view, because 
it is important to consider systems consisting of large but finite number of particles. 
Although there are some works (see \cite{Ec}) 
in this direction, it seems to be hard to resolve.

Thus, the next natural quantity to consider is 
to take an ``average'' of the quantity $\Omega(N,E)$ in the number of particles, $N$. 
So, we define 
\begin{equation}
\label{states1}
\Omega(E)=\sum_{N \geq 1}\Omega(N,E),\ \ E \in \Gamma(P), \quad \Omega(0)=1,
\end{equation}
where we set 
\begin{equation}
\label{discrete}
\Gamma(P)=\{0\} \cup \bigcup_{N \geq 1}{\rm Spec}(P^{N}_{B}) \subset \mb{R}_{\geq 0}. 
\end{equation}
It is easy to see that the set $\Gamma(P)$ is discrete subset of $\mb{R}_{>0}$ without 
finite accumulation points and the sum in \eqref{states1} is finite 
for each fixed $E \in \Gamma (P)$. We call the quantity $\Omega(E)$ the 
{\it number of states of Boson gas}. 
Then our problem is to find asymptotic formulas for $\Omega(E)$ as energy $E$ getting large.

Before stating our result on the asymptotics of $\Omega(E)$, we should mention 
why this kind of problem can be regarded as an analogue of the 
asymptotics of the number of partitions $p(E)$ for positive integer $E$. 
To explain it, consider the generating function $G_{P}$ of the quantity $\Omega(E)$: 
\begin{equation}
\label{gfunc}
G_{P}(\tau)=\sum_{E \in \Gamma(P)}\Omega(E)e^{-E\tau}
\end{equation}
which is called the {\it grand partition function} for Boson gas. Note that the 
usual grand partition function has one more parameter, called chemical potential. 
The chemical potential plays an important and essential role 
in the usual setting-up for analyzing, for example, the Bose-Einstein condensation. 
The main reason for removing the chemical potential is that 
we assumed here that the operator $P$ is strictly positive, and 
so a phenomenon like Bose-Einstein condensation does not occur. 
 
Now it is easy to explain the relation between the number of states $\Omega(E)$ 
of Boson gas defined for $E \in \Gamma(P)$ and 
the number of partitions $p(E)$ for positive integers $E$. 
Namely, we have the following: 
\begin{equation}
\label{infP}
G_{P}(\tau)=\prod_{\ell=1}^{\infty}(1-e^{-\lambda_{\ell}\tau})^{-1}, 
\end{equation}
where $0<\lambda_{1} \leq \lambda_{2} \leq \cdots$ denotes the eigenvalues of 
$P$ counted repeatedly according to their multiplicities (see Lemma \ref{pfunc}).
Thus, if $\lambda_{\ell}=\ell$ for each $\ell$ and the multiplicities
of the eigenvalues are all one, then $\Omega(E)$ is nothing but the 
number of partitions.

We mention here about the motivation 
to address the problem stated in the above. 
One of main motivations comes from the paper \cite{TZ} on 
the various aspects of asymptotics of the multiplicities  $m_{N}(\lambda,\mu)$ of the weights $\mu$ 
in the tensor power $V_{\lambda}^{\otimes N}$ of a fixed irreducible representation $V_{\lambda}$ 
with the dominant weight $\lambda$ of a compact Lie group. 
By the Borel-Weil theorem, 
the irreducible representation $V_{\lambda}$ can be realized as the space $H^{0}(\ocal_{\lambda},L_{\lambda})$ 
of holomorphic sections of a suitable line bundle $L_{\lambda}$ over the coadjoint orbit $\ocal_{\lambda}$ through $\lambda$. 
Then, as a representation, $V_{\lambda}^{\otimes N}$ is naturally identified with $H^{0}(\ocal_{\lambda}^{N},L_{\lambda}^{\boxtimes N})$, 
where $\ocal_{\lambda}^{N}$ is the $N$-times product space of the coadjoint orbit $\ocal_{\lambda}$ 
and $L_{\lambda}^{\boxtimes N}$ is the $N$-th external tensor power of the line bundle $L_{\lambda}$. 
In this point of view, it would be natural to think $\ocal_{\lambda}$ as a classical phase space of a particle, 
and hence $\ocal_{\lambda}^{N}$ would be regarded as a phase space of a system of $N$-particles. 
Then, the multiplicities $m_{N}(\lambda,\mu)$ would be regarded as the number of stationary states 
for the action of the maximal torus. Therefore, the problem on the number of states $\Omega(E)$ of Boson gas 
would be regarded as an analogue of the multiplicities $m_{N}(\lambda,\mu)$ in this context. 
(More precisely, one should consider the multiplicities in the symmetric tensor power ${\rm Sym}^{N}(V_{\lambda})$ 
for a complete analogy of `Boson' gas.)

\subsection{Results}
\label{Intro2}

First of all, we state a theorem which is valid without any assumption on the eigenvalues of $P$. 
In such a general case, it seems to be still hard to obtain the asymptotics of the quantity $\Omega (E)$. 
Thus, we take its average. Define the quantity $D(E)$ by 
\begin{equation}
\label{averaged}
D(E):=\sum_{L \in \Gamma (P) \cap [0,E]}\Omega (L). 
\end{equation}

\begin{maintheo}
\label{WeylAA}
We have 
\begin{equation}
\label{WeylA}
D(E)=e^{B_{n}E^{\frac{n}{n+1}(1+o(1))}},\quad \Gamma (P) \ni E \to \infty, 
\end{equation}
where the constant $B_{n}$ is given by 
\begin{equation}
\label{constB}
B_{n}=(n+1)\left(
\frac{\vol (\Sigma)}{(2\pi n)^{n}}
\zeta (n+1)\Gamma (n)
\right)^{\frac{1}{n+1}}.
\end{equation}
In particular, we have 
\begin{equation}
\label{logaveasympt}
\lim_{\Gamma (P) \ni E \to \infty}
E^{-\frac{n}{n+1}}\log D(E) =B_{n}. 
\end{equation}
\end{maintheo}

In Theorem \ref{WeylAA}, we do not need any assumption on the spectrum 
of $P$ nor classical mechanics. One may think Theorem \ref{WeylAA} as 
an analogue of the usual Weyl asymptotic formula for the operator $P$ itself. 

To state our next results, let us prepare some notation. 
Let 
\begin{equation}
\label{zeta}
Z_{P}(s)=\sum_{\ell \geq 1}\lambda_{\ell}^{-s},\quad s \in \mb{C}
\end{equation}
denote the spectral zeta function. It is well-known (\cite{DG}, \cite{Sh}) that 
the series $Z_{P}(s)$ converges absolutely on the right-half plane $\re (s)>n$, 
and it is meromorphically continued to the whole complex plane. 
The poles of the function $Z_{P}(s)$ are located at 
\begin{equation}
\label{poles}
s=n,n-1,\ldots,1,-1,\ldots,
\end{equation}
and it is holomorphic at $s=0$. We set 
\begin{equation}
\label{residue}
K_{j}=A_{j}\zeta(n-j+1)\Gamma(n-j),\quad j=0,1,\ldots,n-1,\quad 
A_{j}={\rm Res}_{s=n-j}Z_{P}(s),  
\end{equation}
where $\zeta$ and $\Gamma$ denote the Riemann zeta function and the Gamma function, respectively. 
For example, we have $A_{0}=(2\pi )^{-n}\vol (\Sigma)$ with $\Sigma=p^{-1}(1) \subset T^{*}M \setminus 0$ being 
a level surface of the principal symbol $p$ and $\vol(\Sigma)$ being the volume of $\Sigma$ with respect to the Liouville measure. 
Thus, we always have $K_{0}>0$, but, in general, $K_{j}$ could be non-positive (but it is real). 
Let $\varphi_{t}:\Sigma \to \Sigma$ denote the Hamilton flow generated by $p$ (restricted to $\Sigma$). 

\begin{maintheo}
\label{main1}
Assume that the eigenvalues of $P$ are integers. 
Assume also that, for each period $T$ with $0<T \leq \pi$ of a periodic orbit of the 
classical Hamilton flow $\varphi_{t}$, the fixed point set of $\varphi_{T}$ 
is a union of finite number of connected submanifolds and 
each submanifold is clean and of dimension strictly less than $\dim \Sigma=2n-1$. 
Then, we have the following asymptotic formula: 
\begin{equation}
\label{mainF}
\Omega(E)=CE^{\kappa}\exp[x_{E}^{-n}f(x_{E})](1+O(E^{-\kappa_{1}})), 
\end{equation}
where the positive constant $x_{E}$ depending on the energy $E$ satisfies 
\begin{equation}
\label{zeroE}
x_{E}=\left(
\frac{nK_{0}}{E}
\right)^{\frac{1}{n+1}}
(1+O(E^{-\frac{1}{n+1}}))
\end{equation}
$($precise definition of the positive number $x_{E}$ is given in Section {\rm \ref{Proof}}$)$, 
and the function $f(x)$ on the positive real numbers is a polynomial given by 
\begin{equation}
\label{freqF}
f(x)=\sum_{j=0}^{n-1}(n-j+1)K_{j}x^{j}. 
\end{equation}
The constants $C$, $\kappa$ and $\kappa_{1}$ in \eqref{mainF} is given by 
\[
\begin{array}{lll}
C & = & \frac{1}{\det(P)\sqrt{2\pi (n+1)}}(nK_{0})^{\frac{1-2Z_{P}(0)}{2(n+1)}}, \\[10pt]
\kappa & = & \frac{Z_{P}(0)-1-n/2}{n+1}, 
\end{array}
\]
where $\det (P)$ is the zeta regularized determinant of $P$ and the constant $\kappa_{1}$ is 
any number satisfying $0<\kappa_{1}<\frac{1}{n+1}$ when $n\geq2$, and $0<\kappa_{1}<1/4$ when $n=1$. 
\end{maintheo}

Recall that a submanifold $Z$ in $\Sigma$ of fixed points of 
a diffeomorphism $\varphi:\Sigma \to \Sigma$ 
is said to be clean (\cite{DG}) if, for each $z \in Z$, the set of fixed points 
of the differential $d\varphi_{z}:T_{z}\Sigma \to T_{z}\Sigma$ equals the tangent space 
of $Z$ at $z$. A reason why we assume a cleanness condition for the classical 
Hamilton flow in Theorem \ref{main1} is that we need to use a theorem in \cite{DG} on 
singularities of the trace of the wave operator $e^{-itP}$. See Sections \ref{Intro4}, \ref{Properties2}. 

In Theorem \ref{main1}, the assumption that the eigenvalues of $P$ are integers implies 
that the classical Hamilton flow $\varphi_{t}$ on $T^{*}M$ is periodic with common period $2\pi$. 
This follows from Egorov's theorem. 
But, there might be other period less than $2\pi$. The assumption that the fixed point set 
of $\varphi_{T}$ for a period $T$ of a periodic orbit is of dimension strictly less than $2n-1$ 
implies that $2\pi$ is the least common period of the flow. There is a technical reason 
why common periods less than $2\pi$ are excluded. See Section \ref{Properties2}. 

One of examples of the operator $P$ satisfying the assumption in Theorem \ref{main1} 
is the operator $P=\sqrt{\Delta_{{\rm can}}+(n-1)^{2}/4}+(n-1)/2$ with 
the standard Laplacian $\Delta_{{\rm can}}$ on the $n$-sphere. 

Note also that, according to a result of Colin de Verdier\`{e} (\cite{CdV}), if the classical 
flow is periodic with the least common period $2\pi$, and the square $P^2$ of 
the operator $P$ is a differential operator, then there is a constant $\alpha$ and 
a pseudo-differential operator $Q$ of order $-1$ such that 
the eigenvalues of the operator $P+Q-\alpha/4$ are integers. 
Therefore, we can apply Theorem \ref{main1} for these systems if 
the cleanness assumption is satisfied.  

Compared with the usual Weyl asymptotic formula for one particle, 
the asymptotic formula for the number of states $\Omega(E)$ of Boson gas 
is rather complicated. 
In particular, all the positive poles affects the exponent of the asymptotics 
whereas only the largest pole affects in the usual Weyl asymptotic formula. 
As is shown in Section \ref{Proof}, the positive constant $x_{E}$ satisfies 
the estimate \eqref{zeroE}. 
Therefore, if one takes the logarithm and then the limit to reduce information, 
one obtain the following. 
\begin{maincor}
\label{limL}
Under the same assumption as in Theorem {\rm \ref{main1}}, we have 
\begin{equation}
\label{Knopp}
\lim_{\Gamma (P) \ni E \to \infty}E^{-\frac{n}{n+1}}\log \Omega(E)
=B_{n},  
\end{equation}
where the constant $B_{n}$ is given in $\eqref{constB}$. 
\end{maincor}

Corollary \ref{limL} is basically obtained by Knopp (\cite{Kn}), where the 
sequence $\{\lambda_{n}\}$ of positive integers is assumed to be strictly increasing, 
that is, all the eigenvalues of the operator $P$ are of multiplicity one. 
Thus, one can think Corollary \ref{limL} as a tiny generalization of Knopp's result. 
One might expect that the limit formula \eqref{Knopp} might hold for general system 
without the assumption that the eigenvalues of $P$ are integers. 
However, in general, since $\Omega (E)$ does not have monotonicity in $E$, it might be hard to 
obtain it. 

In Theorem \ref{main1}, the positive constant $x_{E}$ is defined as a zero of a 
polynomial determined by the residues of the positive poles of the spectral zeta function. 
Thus, it is not so transparent. However, when $n=\dim M=2$, it is replaced by a 
more concrete constant $y_{E}$ as follows. 

\begin{maintheo}
\label{main2}
Let $n=\dim M =2$. Then, under the same assumption in Theorem {\rm \ref{main1}}, 
we have the following. 
\[
\Omega(E)=CE^{\frac{Z_{P}(0)-2}{3}}\exp
\left[
\frac{3K_{0}}{y_{E}^{2}}+\frac{K_{1}}{y_{E}}-\frac{K_{1}^{2}}{12K_{0}}
\right] (1+O(E^{-\kappa_{1}}))
\]
for any $0<\kappa_{1}<1/3$, where the constants $C$ and $y_{E}$ is given by 
\[
y_{E} = \left(
\frac{\zeta(3)\vol(\Sigma)}{2\pi^{2}E}
\right)^{1/3},
\qquad
C = \frac{1}{\det (P)\sqrt{6\pi}}
\left(
\frac{\zeta(3)\vol (\Sigma)}{2\pi^{2}}
\right)^{(1-2Z_{P}(0))/6}. 
\]
\end{maintheo}

\subsection{Meinardus' theorem}

Among results on asymptotics of the number of partitions, 
a theorem due to Meinardus \cite{Me} gives a natural 
generalization of the above formula \eqref{hra} of Hardy-Ramanujan. 
Let us state a result of Meinardus here. 
Let $\{a_{\ell}\}_{\ell=1}^{\infty}$ be a sequence of non-negative numbers (not necessarily integers) 
such that the infinite product 
\[
G(\tau)=\prod_{\ell \geq 1}(1-e^{-\ell\tau})^{-a_{\ell}}
\]
converges locally uniformly for $\re (\tau) >0$. Then the infinite product $G(\tau)$ defines 
a holomorphic function on the right-half plane. Thus we can consider the Taylor expansion 
of the function $G(z)=G(\tau)$ ($z=e^{-\tau}$) at $z=0$: 
\[
G(z)=1+\sum_{E=1}^{\infty}r(E)z^{E},\quad z=e^{-\tau},\quad \re (\tau)>0. 
\]
The Meinardus theorem is an asymptotic formula for the coefficients $r(E)$ as 
$E \to \infty$. When $a_{n}=1$ for all $n$, usual number of partitions $p(E)$ coincides with 
the quantity $r(E)$. 
Associated to the sequence $\{a_{n}\}$ is the Dirichlet series
\begin{equation}
L(s)=\sum_{\ell=1}^{\infty}\frac{a_{\ell}}{\ell^{s}},\quad s \in \mb{C}. 
\end{equation}
To state the Meinardus theorem, we need the following assumptions: 

\begin{enumerate}
\item There exists a positive number $\alpha$ such that $L(s)$ converges for $\re (s) >\alpha$. 
\item There exists a number $C_{0}$ such that $0 < C_{0} < 1$ and the Dirichlet series $L(s)$ is continued 
meromorphically on the set $\re (s) \geq -C_{0}$. 
\item The function $L(s)$ has just one pole at $s=\alpha$ in $\re (s) \geq -C_{0}$, 
and it is simple with residue $A={\rm Res}_{s=\alpha}L(s)>0$. 
\item There exists a number $R>0$ such that $|L(s)| =O(|\im (s)|^{R})$ for $\re(s) \geq -C_{0}$. 
\item[(H)] Define the function $\theta (\tau)$ ($\re(\tau) >0$) by
\begin{equation}
\theta (\tau)=\sum_{\ell=1}^{\infty}a_{\ell}e^{-\ell\tau}. 
\end{equation}
Then, for every sufficiently small $\ep >0$, there exists a constant $C>0$ such that 
for every $\tau=x+iy$ with $|\arg (\tau)| >\pi/4$, $|y| \leq \pi$, we have 
\begin{equation}
\label{theta0}
\re (\theta (x+iy)) -\theta (x) \leq -Cx^{-\ep}. 
\end{equation}
\end{enumerate}
Note that, by the assumption (1), $L(s)$ is holomorphic on $\re(s) >\alpha$. 
The assumption (4) is to guarantee the convergence of an integral. 
Then, the Meinardus theorem is stated as follows. 

\begin{maintheo}[Meinardus (\cite{Me})]
\label{meinardusT}
Assume that the sequence $\{a_{n}\}_{n=1}^{\infty}$ of non-negative numbers satisfies the above 
conditions {\rm (1) -- (4)} and {\rm (H)}. Then, we have the following asymptotic formula for $r(E)$. 
\[
r(E)=CE^{\kappa}\exp
\left[
\left(
\frac{\alpha +1}{\alpha}
\right)
E^{\frac{\alpha}{\alpha +1}}
(A\Gamma(\alpha +1)\zeta(\alpha +1))^{\frac{1}{\alpha +1}}
\right]
\times (1+O(E^{-\kappa_{1}})),  
\]
where $A>0$ is the residue of $L(s)$ at $s=\alpha$, 
and the constants $C$, $\kappa$ and $\kappa_{1}$ are given by 
\[
\begin{array}{lll}
C & = & \frac{e^{L'(0)}}{\sqrt{2\pi (\alpha +1)}}
(A\Gamma(\alpha +1)\zeta(\alpha +1))^{\frac{1-2L(0)}{2(\alpha +1)}} \\
\kappa & = & \frac{L(0)-1-\alpha/2}{1+\alpha} \\
\kappa_{1} & = & \frac{\alpha}{\alpha +1}\min 
\left\{
\frac{C_{0}}{\alpha}-\frac{\delta}{4},\ \frac{1}{2}-\delta
\right\}\\
 & & 0 < \delta < \frac{1}{2}
\end{array}
\]
\end{maintheo}

Meinardus' method of obtaining the above theorem 
is based on the saddle point method, and thus it would be natural to expect that 
Meinardus' method would work well for more general problem. 
Roughly speaking, 
Meinardus approximates the function $\log G(\tau)$, which plays a role of the phase function, 
by a very simple function and then estimates the error terms 
by using the assumption (H).

\subsection{Remarks and comments}
\label{Intro4}

\noindent{(1)} It would be clear from the statement of Theorem \ref{main1} 
that the spectral zeta function $Z_{P}(s)$ plays the same role 
as the Dirichlet series $L(s)$ in Theorem \ref{meinardusT}. 
Furthermore, the eigenvalues of $P$ are assumed to be integers. 
Therefore, the situation in Theorem \ref{main1} is, at least formally, contained in that of the Meinardus theorem 
with $\{a_{n}\}$ being the multiplicities of the eigenvalues. 
However, there are something to be careful. One of them is the fact that 
the spectral zeta function could have many positive real poles, although 
it is assumed, in Theorem \ref{meinardusT}, that $L(s)$ has only one positive pole. 
Another is the condition (H) which is assumed in Theorem \ref{meinardusT}. 
In Theorem \ref{main1}, there are no assumption similar to the condition (H). 
In our case, the condition (H) is a condition on the theta function (`heat trace')
\begin{equation}
\label{theta}
\theta (\tau)=\sum_{\ell =1}^{\infty}e^{-\lambda_{\ell}\tau},\quad \re (\tau) >0. 
\end{equation}
Then, the cleanness assumption for the classical mechanics in Theorem \ref{main1} 
makes us to be able to prove the condition (H) for the function $\theta$ (see Section \ref{Properties}). 

The assumptions in Theorem \ref{main1} exclude 
the existence of the common period less than $2\pi$. 
This is necessary for the validity of the condition (H) for the function $\theta$ defined above. 
In fact, if all the eigenvalues are even integers, then we have $\theta(x+i\pi)=\theta(x)$ 
for $x>0$. This shows that the condition (H) does not hold for the system having eigenvalues all 
of them are even integers. By Helton-Guillemin theorem (\cite{Gu2}), 
if $\pi$ is a common period of the classical Hamilton flow, then (assuming that $P$ being the 
square root of the Laplacian), the set of cluster points of differences of eigenvalues of $P$, 
$\{\lambda_{i}-\lambda_{j}\}$, is just $2\mb{Z}$. Thus, a perturbation of $P$ could have eigenvalues all 
of them are even integers. Thus, if $\pi$ is a common period of the classical flow, 
then, one might need to impose some additional assumption on the eigenvalues of $P$ to 
ensure the validity of the condition (H).

\vspace{10pt}

\noindent{(2)}  As mentioned, we assumed that the operator $P$ is positive definite, 
and hence a phenomenon like Bose-Einstein condensation does not occur. 
However, when the operator $P$ has the zero eigenvalue, 
or when we consider a family of positive definite operators $P_{t}$ 
with the least eigenvalues $\lambda_{1}(t)$ which tends  to zero as $t \to \infty$, 
then the number of particles 
in the zero energy state, or small energy states, could be large. 
This phenomenon is Bose-Einstein condensation, 
and in such a case, we need to use the full grand partition function defined as 
\[
\Xi(\tau,\mu)=
\prod_{\ell=1}^{\infty}
(1-e^{\mu-\lambda_{\ell}\tau})^{-1},
\]
where the parameter $\mu$, which is called the chemical potential, should be negative. 
In this case, the quantity $\Omega(E)$ should be replaced by 
\[
\Omega(E,\mu)=\sum_{N \geq 1}\Omega(N,E)e^{N\mu}. 
\]
Thus, it might be interesting to find joint asymptotics of $\Omega(E,\mu)$ as $E \to \infty$ and $\mu \to 0$. 

\vspace{10pt}

\noindent{(3)} It would be hard to describe the 
asymptotics of $\Omega(E)$ for general system. 
This is because, in general, there are no integral representation 
for the quantity $\Omega(E)$. Although there are some representation for $\Omega(E)$ in terms of 
a limit of integrals (see Lemma \ref{intL}), we need to estimate the error term. 
In the error term, such a quantity like the difference of eigenvalues in $\Gamma(P)$ appears. 
In general, estimating the level spacing is quite hard. 
However, when the classical system is periodic, like a Zolll surface, one might 
be able to use the results in \cite{DG}, \cite{Gu}, \cite{UZ} to find a similar asymptotic formula 
for $\Omega(E)$ in such a case. 

\vspace{10pt}

\noindent{{\bf Acknowledgments}} The author would like to thank to professors T.~Matsui and S.~Zelditch 
for their helpful comments. The author also would like to thank to the referee who 
pointed out some mistakes in the earlier version of the paper and gave valuable suggestions to the author.

\section{Partition function for Boson Gas}
\label{Partition}

In this section, we give precise definition and some properties 
of the number of states, $\Omega(E)$, of Boson gas and 
its partition function, which are, at least formally, quite well known in 
physics textbooks on quantum statistics (\cite{Fu}, \cite{Kh}). 

\subsection{Occupation number representation}
\label{Partition1}

As in Introduction, let $P$ be an elliptic, strictly positive 
pseudo-differential operator of order $1$ on a compact 
connected manifold $M$ of dimension $n$, and consider, 
for each positive integer $N$, the operator $P^{N}$ on $L^{2}(M^{N})$ 
where $M^{N}=M \times \cdots \times M$ ($N$ times) defined in \eqref{gop}. 
We define the action of the symmetric group $\mathfrak{S}_{N}$ on $L^{2}(M^{N})$ by 
\[
\begin{gathered}
(U_{\sigma} F)(x_{1},\ldots,x_{N})=F(x_{\sigma(1)},\ldots,x_{\sigma(N)}),\\
\sigma \in \mathfrak{S}_{N},\ F \in L^{2}(M^{N}), \ (x_{1},\ldots,x_{N}) \in M^{N}. 
\end{gathered}
\]
This action of $\mathfrak{S}_{N}$ on $L^{2}(M^{N})$ commutes with 
the operator $P^{N}$. 
Let $\hcal_{S}^{N}$ and $\hcal_{A}^{N}$ be the space of 
symmetric and anti-symmetric functions on $M^{N}$, respectively: 
\[
\begin{gathered}
\hcal_{S}^{N}=\{F \in L^{2}(M^{N})\,;\,U_{\sigma} F=F,\ \sigma \in \mathfrak{S}_{N}\},
\\
\hcal_{A}^{N}=\{F \in L^{2}(M^{N})\,;\,U_{\sigma} F=\sgn (\sigma)F,\ 
\sigma \in \mathfrak{S}_{N}\}.
\end{gathered}
\]
The spaces $\hcal_{S}^{N}$ and $\hcal_{A}^{N}$ are 
the state space for the system of $N$-Boson and $N$-Fermion, respectively. 
In this paper, we consider the system of Boson. 
We denote the restriction of $P^{N}$ on $\hcal_{S}^{N}$ by $P^{N}_{B}$.

Let $0 < \lambda_{1} \leq \lambda_{2} \leq \cdots \uparrow \infty$ be 
the sequence of eigenvalues of $P$. We choose an orthonormal 
basis $\{\varphi_{j}\}_{j \geq 1}$ of $L^{2}(M)$ consisting of eigenfunctions of $P$:  
$P\varphi_{j}=\lambda_{j}\varphi_{j}$, $j=1,2,\ldots$.
For any multi-index $\alpha =(\alpha_{1},\ldots,\alpha_{N}) \in \mb{Z}_{>0}^{N}$, 
we define a function $\Phi_{\alpha} \in C^{\infty}(M^{N})$ by 
\[
\Phi_{\alpha}=\varphi_{\alpha_{1}}\otimes \cdots \otimes 
\varphi_{\alpha_{N}}. 
\]
Then, the function $\Phi_{\alpha}$ is an eigenfunction of $P^{N}$, and we have 
\begin{equation}
\label{eigenF1}
P^{N}\Phi_{\alpha}=E\Phi_{\alpha},\quad 
E=\sum_{i=1}^{N}\lambda_{\alpha_{i}}. 
\end{equation}
From this, the spectrum of $P^{N}$ is easily described. 
But, what we need to know is the spectrum of the Bose Hamiltonian $P^{N}_{B}$. 
Let $\pi^{N}_{S}$ denote the orthonormal projection onto $\hcal^{N}_{S}$: 
\[
\pi^{N}_{S}=\frac{1}{N!}\sum_{\sigma \in \mathfrak{S}_{N}}U_{\sigma}.
\]
Then, clearly the state space for Boson $\hcal^{N}_{S}$ spanned by 
the functions $\pi^{N}_{S}\Phi_{\alpha}$, $\alpha \in \mb{Z}_{>0}^{N}$. 
The symmetric group $\mathfrak{S}_{N}$ acts on $\mb{Z}_{>0}^{N}$ by 
\[
\sigma \alpha =(\alpha_{\sigma^{-1}(1)},\ldots,\alpha_{\sigma^{-1}(N)}),\quad 
\sigma \in \mathfrak{S}_{N},\ \alpha =(\alpha_{1},\ldots,\alpha_{N}) 
\in \mb{Z}_{>0}^{N}.
\]
\begin{lem}
\label{pre11}
\begin{enumerate}
\item $\pi_{S}^{N}\Phi_{\alpha}=\pi_{S}^{N}\Phi_{\beta}$ if and only if 
$\beta$ lies in the orbit of $\mathfrak{S}_{N}$ through $\alpha$. 
\item Let $[\alpha]$ denote the orbit of $\mathfrak{S}_{N}$ 
through $\alpha \in \mb{Z}_{>0}^{N}$, and set 
$\Phi_{[\alpha]}=\pi_{S}^{N}\Phi_{\alpha}/\|\pi_{S}^{N}\Phi_{\alpha}\|$. 
Then, $\{\Phi_{[\alpha}]\}_{[\alpha] \in \mb{Z}_{>0}^{N}/\mathfrak{S}_{N}}$ form 
an orthonormal basis for $\hcal^{N}_{S}$. 
\end{enumerate}
\end{lem}
\begin{proof}
First, note that we have $U_{\sigma}\Phi_{\alpha}=\Phi_{\sigma \alpha}$ for 
any $\sigma$ and $\alpha$. 
Now, suppose that $\beta =\mu \alpha$ for some $\mu \in \mathfrak{S}_{N}$. 
Then, we have 
\[
\pi_{S}^{N}\Phi_{\beta}=\frac{1}{N!}\sum_{\sigma \in \mathfrak{S}_{N}}
U_{\sigma}\Phi_{\mu \alpha}
=\frac{1}{N!}\sum_{\sigma} \Phi_{\sigma \mu \alpha}=
\frac{1}{N!}\sum_{\sigma}\Phi_{\sigma \alpha}=\pi_{S}^{N}\Phi_{\alpha}.  
\]
Conversely, let $\pi_{S}^{N}\Phi_{\alpha}=\pi_{S}^{N}\Phi_{\beta}$. 
Since $\{\Phi_{\alpha}\}_{\alpha \in \mb{Z}_{>0}^{N}}$ is an orthonormal basis 
for $L^{2}(M^{N})$, we have 
\[
\ispa{\pi_{S}^{N}\Phi_{\alpha},\,\Phi_{\alpha}}
=\frac{1}{N!}\sum_{\sigma}\ispa{\Phi_{\sigma \alpha},\,\Phi_{\alpha}}
=\frac{1}{N!}\sharp{\rm Stab}(\alpha), 
\]
where ${\rm Stab}(\alpha)$ denotes the stabilizer of $\alpha$ in $\mathfrak{S}_{N}$. 
In particular, $\pi_{S}^{N}\Phi_{\alpha}$ is not zero, because $\sharp{\rm Stab}(\alpha)$ is greater than or equal to one. 
Thus, by the assumption that $\pi_{S}^{N}\Phi_{\beta}=\pi_{S}^{N}\Phi_{\alpha}$, 
we have 
\[
\frac{1}{N!}\sum_{\sigma}\ispa{\Phi_{\sigma \beta},\,\Phi_{\alpha}}
=\ispa{\pi_{S}^{N}\Phi_{\beta},\,\Phi_{\alpha}}=\frac{1}{N!}\sharp{\rm Stab}(\alpha), 
\]
and we have $\sharp{\rm Stab}(\alpha) \geq 1$. 
Thus, for some $\sigma \in \mathfrak{S}_{N}$ in the sum of the above, 
$\ispa{\Phi_{\sigma \beta},\,\Phi_{\alpha}}$ is not zero, 
and hence $\alpha=\sigma \beta$. 
This shows (1). (2) follows easily from (1). 
\end{proof}

Since we clearly have $P^{N}_{S}\Phi_{[\alpha]}=E\Phi_{[\alpha]}$ 
with $E$ given in \eqref{eigenF1}, the spectrum ${\rm Spec}(P_{B}^{N})$ 
of $P_{B}^{N}$ is descrete and coincides with ${\rm Spec}(P^{N})$. 
Thus, what is important is to understand the multiplicities of eigenvalues of $P_{B}^{N}$.

For any $E \in {\rm Spec}(P^{N}_{B})$, we set 
\[
C_{N}(E)=\{\alpha \in \mb{Z}_{>0}^{N}\,;\,
E=\sum_{j=1}^{N}\lambda_{\alpha_{j}}\}. 
\]
We then define 
\[
\Omega(N,E):=\sharp[C_{N}(E)/\mathfrak{S}_{N}]. 
\]
Then, by Lemma \ref{pre11} (2), the quantity $\Omega (N,E)$ equals 
the multiplicity of the eigenvalue $E$ of the Bose Hamiltonian $P_{B}^{N}$. 

\begin{lem}
\label{occupy}
For any $E \in {\rm Spec}(P_{B}^{N})$, we have 
\begin{equation}
\label{orep}
\Omega(N,E)=\sharp
\left\{
(n_{k})_{k=1}^{\infty} \in \mb{Z}_{\geq 0}^{\infty}\,;\,
E=\sum_{k \geq 1}n_{k}\lambda_{k},\ N=\sum_{k \geq 1}n_{k}
\right\}
\end{equation}
\end{lem}
\begin{proof}
For each $\alpha \in C_{N}(E)$, 
we define non-negative integers $n_{k}(\alpha)$ by 
\[
n_{k}(\alpha)=\sharp\{j \,;\,\alpha_{j}=k\}. 
\]
Then, we have 
\[
E=\sum_{j=1}^{N}\lambda_{\alpha_{j}}=\sum_{k=1}^{\infty}n_{k}(\alpha)\lambda_{k}, 
\quad \sum_{k=1}^{\infty}n_{k}(\alpha)=N.  
\]
It is easy to see that, for $\alpha,\beta \in \mb{Z}_{>0}^{N}$, 
$n_{k}(\alpha)=n_{k}(\beta)$ for every $k$ if and only if $\beta$ 
is in the $\mathfrak{S}_{N}$-orbit through $\alpha$, and from this fact, 
we conclude the assertion. 
\end{proof}

The expression of $\Omega(N,E)$ in \eqref{orep} is called the {\it occupation number representation}. 
By Lemma \ref{occupy}, the sum in the definition \eqref{states1} of $\Omega(E)$ 
is a finite sum for each fixed $E$, if the eigenvalues of $P$ are positive. 
The quantity $\Omega(E)$, with setting $\Omega(0)=1$, is thus a function on the set 
$\Gamma(P)$ defined in \eqref{discrete}. 
For each positive number $A$, the set $\Gamma(P) \cap [0,A]$ is finite because 
if we write $E \in \Gamma (P) \cap [0,A]$ as $E=\sum_{k}n_{k}\lambda_{k}$ with $N=\sum_{k}n_{k}$ as in \eqref{orep}, 
we clearly have $\lambda_{1}N \leq E \leq A$. Therefore, $\Gamma (P)$ is a discrete subset in $\mb{R}_{\geq 0}$ without 
finite accumulation points.

\subsection{Partition function}
\label{Partition2}

The partition function for the number of states $\Omega(E)$ of Boson gas is defined, at least formally, by the 
equation \eqref{gfunc}. 
The function $G_{P}(\tau)$ always converges and has infinite product 
representation if $P$ is strictly positive definite as follows. 

\begin{lem}
\label{pfunc}
The partition function $G_{P}(\tau)$ defined in {\rm \eqref{gfunc}} converges absolutely and 
locally uniformly for $\re(\tau) >0$, and on this region, $G_{P}(\tau)$ has the infinite product representation {\rm \eqref{infP}}.
\end{lem}
\begin{proof}
Set $u_{\ell}(\tau) =(1-e^{-\lambda_{\ell}\tau})^{-1}-1$. Then we have 
\[
|u_{\ell}(\tau)| \leq (1-e^{-\lambda_{1}x})^{-1} e^{-\lambda_{\ell}x}, 
\]
where we set $\tau =x+iy$. Since the heat trace $\theta (x)=\sum_{\ell \geq 1}e^{-\lambda_{\ell}x}$ converges 
locally uniformly for $x>0$ (see, for example, \cite{Gi}), 
the sum $\sum u_{\ell}(\tau)$ converges absolutely and  locally uniformly for $\re (\tau) >0$. 
Therefore, the infinite product in \eqref{infP} converges absolutely and defines a holomorphic function on the 
region $\re (\tau) >0$. 
Now, we compute the partial product of the 
infinite product in \eqref{infP}, by using the geometric series expansion for each term, 
as follows. 
\[
\prod_{\ell=1}^{K}(1-e^{-\lambda_{\ell}\tau})^{-1}
=\sum_{E \in \Gamma (P)}\Omega_{K}(E)e^{-E\tau}, 
\]
where we set
\[
\Omega_{K}(E)=\sharp
\left\{
(n_{k})_{k=1}^{K} \in \mb{Z}_{\geq 0}^{K}\,;\,
\sum_{k=1}^{K}n_{k}\lambda_{k}=E
\right\}. 
\]
For each $L >0$, there exists $K_{L} >0$ such that $L < \lambda_{K_{L}+1}$. 
Then, it is easy to see that $\Omega(E)=\Omega_{K_{L}}(E)$ for each $E \in \Gamma(P)$, $E <L$. Therefore, 
for positive $x$, we have 
\[
\prod_{\ell=1}^{\infty}(1-e^{-\lambda_{\ell}x})^{-1} \geq 
\sum_{E \in \Gamma(P),\,E<L}\Omega(E)e^{-Ex}
\]
for each $L>0$. Thus the partition function $G_{P}(\tau)$, defined in \eqref{gfunc}, 
converges absolutely and locally uniformly for $\re(\tau)>0$, and satisfies 
\[
G_{P}(x) \leq \prod_{\ell=1}^{\infty}(1-e^{-\lambda_{\ell}x})^{-1}
\]
for positive number $x$. Conversely, it is clear that 
for each $K>0$ and $E \in \Gamma(P)$, we have 
\[
\prod_{\ell=1}^{K}(1-e^{-\lambda_{\ell}x})^{-1}
\leq G_{P}(x), 
\]
and which, with the analytic continuation, proves the lemma. 
\end{proof}

If the eigenvalues of $P$ are integers, then we have 
\begin{equation}
\label{int2}
\Omega(E) =\frac{1}{2\pi}\int_{|y| \leq \pi}
e^{E(x+iy)}G_{P}(x+iy)\,dy
\end{equation}
for any positive number $x$. This formula is the starting point of our analysis. 
In general, the set $\Gamma (P)$ is not contained in the set of integers, and hence 
there are no such an integral representation for $\Omega(E)$. 
However, even if $\Gamma (P) \not \subset \mb{Z}$, we can approximate $\Omega(E)$ 
by integrals of the form \eqref{int2}. We shall mention about this formula, although 
this is not necessary for later sections. 
\begin{lem}
\label{intL}
For any $x>0$ and $T>0$, we have the following. 
\begin{equation}
\label{asint}
\left|
\Omega(E) -\frac{1}{2T}
\int_{-T}^{T}
e^{E(x+iy)}G_{P}(x+iy)\,dy
\right|
\leq \frac{G_{P}(x)e^{Ex}}{T}
\left(
\min_{E' \in \Gamma (P),E' \neq E}|E-E'|
\right)^{-1}
\end{equation}
\end{lem}
\begin{proof}
First of all, we note that the minimum in the right hand side in \eqref{asint} exists 
because $\Gamma (P)$ is a discrete subset of $\mb{R}_{>0}$ without finite accumulation points. 
Substituting the series expansion for $G_{P}(x+iy)$ for the integral in the left-hand side 
of \eqref{asint}, we have 
\[
\frac{1}{2T}\int_{-T}^{T}e^{E(x+iy)}G_{P}(x+iy)\,dy -\Omega(E)
=\sum_{E' \in \Gamma (P),E' \neq E}
\Omega(E')\frac{\sin (T(E-E'))}{T(E-E')}
e^{(E-E')x}, 
\]
and this proves \eqref{asint}. 
\end{proof}

\section{Properties of the partition function}
\label{Properties}

In this section, we review some facts about the spectral zeta functions 
defined in \eqref{zeta}, and then, by using these facts, 
we derive an approximation formula for the partition function, which is a tiny generalization 
of a lemma in \cite{Me}, \cite{An}. 
We prove the condition (H) under the periodicity assumtion for the classical 
dynamical system.

\subsection{Approximation of the partition function}
\label{Properties1}

It is well-known (\cite{DG}, \cite{Sh}) that the spectral zeta function $Z_{P}(s)$ for 
the first order elliptic positive pseudo-differential operator $P$ has the following properties. 
\begin{enumerate}
\item The series $Z_{P}(s)$ converges absolutely for $\re(s) >n$, and it is continued meromorphically 
on the whole complex plane. 
\item All of the poles of the function $Z_{P}(s)$ are simple and they are given by \eqref{poles}. 
In particular, $Z_{P}(s)$ is holomorphic around the origin. 
\item For each $c \in \mb{R}$, the spectral zeta function $Z_{p}(s)$ has a polynomial growth 
in $\im (s)$ locally uniformly in $\re (s) \geq c$ (except on neighborhoods of poles). 
\end{enumerate}

By using these facts about the spectral zeta function $Z_{P}(s)$ we have 
the following. 

\begin{lem}
\label{formula1}
Let $K_{j}$ be the real numbers defined in {\rm \eqref{residue}}. Then, for $\tau=x+iy \in \mb{C}$ with 
$0<x$, $|y| \leq x$, the logarithm of the partition function $\log G_{P}(\tau)$ has the following expression: 
\begin{equation}
\label{logP}
\log G_{P}(\tau)=\sum_{j=0}^{n-1}K_{j}\tau^{-(n-j)}
-Z_{P}(0)\log \tau +Z_{P}'(0)+J(\tau),
\end{equation}
where the term $J(\tau)$ is given by 
\begin{equation}
\label{error0}
J(\tau)=\frac{1}{2\pi i}\int_{-C_{0}-i\infty}^{-C_{0}+i\infty}
\tau^{-s}Z_{P}(s)\zeta(s+1)\Gamma (s)\,ds
\end{equation}
for any $0 <C_{0} <1$. Furthermore, for $\tau=x+iy$ with $0<x$, $|y| \leq x$, we have 
$|J(\tau)|=O(x^{C_{0}})$. 
The branch of the logarithm in \eqref{logP} is described in the proof.
\end{lem}
\begin{proof}
Proof is the same as that in \cite{Me}, \cite{An}, but we review it for completeness. 
According to the infinite product representation \eqref{infP} of the partition function $G_{P}(\tau)$, 
we can define the branch of the logarithm $\log G_{P}(\tau)$ by 
\[
\log G_{P}(\tau)=-\sum_{\ell=1}^{\infty}\log (1-e^{-\lambda_{\ell}\tau}), 
\]
where the logarithm in the right-hand side is its principal branch. 
For $\re (\tau)>0$ and $\sigma_{0}>0$, we have 
\[
e^{-\tau}=\frac{1}{2\pi i}\int_{\sigma_{0} -i\infty}^{\sigma_{0} +i\infty}
\tau^{-s}\Gamma (s)\,ds.   
\]
By using the above and the Taylor expansion of the logarithm, 
\begin{equation}
\label{logT}
-\log (1-e^{-\lambda_{\ell}\tau})=\sum_{k \geq 1}
\frac{1}{k}e^{-\lambda_{\ell}k\tau} \quad (\re (\tau) >0), 
\end{equation}
we have 
\[
-\log (1-e^{-\lambda_{\ell}\tau})
=\frac{1}{2\pi i}
\int_{\sigma_{0}-i\infty}^{\sigma_{0}+i\infty}
\lambda_{\ell}^{-s}\tau^{-s}\zeta (s+1)\Gamma (s)\,ds 
\]
for $\re (\tau) >0$ and $\sigma_{0} >0$. Taking $\sigma_{0} >n$ and summing over $\ell \geq 1$, 
we have 
\[
\log G_{P}(\tau)=\frac{1}{2\pi i}
\int_{\sigma_{0} -i\infty}^{\sigma_{0}+i\infty}
\tau^{-s}Z_{P}(s)\zeta (s+1)\Gamma (s)\,ds. 
\]
Take $0 < C_{0} <1$ and set $s=\sigma +it$ with $-C_{0} \leq \sigma \leq \sigma_{0}$ and $t \in \mb{R}$. 
Then, for $\tau=x+iy$ with $0<x$, $|y| \leq x$ and $s=\sigma +it$ with the above region, we have 
\begin{equation}
\label{auxest1}
|\tau^{-s}| \leq |\tau|^{-\sigma}e^{\pi |t|/4}. 
\end{equation}
It is well-known that $\zeta (1+s)$ is of order $O(|t|^{c_{1}})$ with a positive constant $c_{1}$. 
By Corollary 2.2 in \cite{DG}, we have $|Z_{P}(s)| =O(|t|^{c_{2}})$ with a positive constant $c_{2}$. 
It is also well-known that we have $|\Gamma (s)| =O(|t|^{c_{3}}e^{-\pi |t|/2})$ with a positive constant $c_{3}$. 
Thus, we can use the residue formula to shift the integral contour to $-C_{0}+it$, $t \in \mb{R}$. 
The function $f(s):=\tau^{s}Z_{P}(s)\zeta (s+1)\Gamma (s)$ has 
simple poles $s$ with $-C_{0} \leq \re (s) \leq \sigma$ at $s=n-j$ with $j=0,\ldots , n-1$ 
and a pole of order two at $s=0$. 
A direct computation shows that the residue of $f(s)$ at $s=n-j$ $(j=0,\ldots,n-1)$ is $K_{j}\tau^{-(n-j)}$ 
and the residue of $f(s)$ at $s=0$ is $-Z_{P}(0)\log \tau +Z_{P}'(s)$. 
Thus, the residue formula gives the expression \eqref{logP} with the function $J(\tau)$ given by \eqref{error0} 
for $\tau=x+iy$ with $x>0$ and $|y| \leq x$. Note that the integral $J(\tau)$ is absolutely convergent 
for this region and is of order $O(x^{C_{0}})$ by \eqref{auxest1}, which completes the proof. 
\end{proof}

\subsection{The condition (H)}
\label{Properties2}

Our next purpose is to prove the condition (H) on the function $\theta (\tau)$ defined by \eqref{theta}
under a cleanness assumption of the classical Hamilton flow $\varphi_{t}$. 
First, we begin with a heuristic argument by using the heat trace asymptotics
\begin{equation}
\label{heat1}
\theta(\tau) \sim c_{0}\tau^{-n}
\end{equation}
whose validity is quite well-known at least for positive $\tau$. 
If \eqref{heat1} would hold for $\tau \in \mb{C}$, $\re (\tau) >0$ as $\tau \to 0$, 
the function $\re \theta (\tau) -\theta (\re (\tau))$ in the condition (H) would 
be replaced by more simple function $c_{0}(\re (\tau^{-n})-\re (\tau)^{-n})$. 
Then, the condition (H) would be proved, at least for small $|y|$, by the inequality 
\begin{equation}
\label{aux11}
\re ((x+iy)^{-n}) -x^{-n} \leq -\frac{1}{2}x^{-n}
\end{equation}
for $0<x \leq |y|$. 
However, when $\re (\tau)>0$ and $\tau \to 0$, we face singularities of 
the distribution $\theta (it)$, $t \in \mb{R}$, on the imaginary axis, and hence we need to 
take care about such singularities to prove the condition (H).

\begin{prop}
\label{condH}
Assume that, for each period $T$ with $0< T \leq \pi$ of a periodic orbit of the flow $\varphi_{t}$, the fixed point set 
of $\varphi_{T}:\Sigma \to \Sigma$ is a union of finite number of connected submanifolds $Z_{1},\ldots,Z_{r}$ and 
each $Z_{j}$ is clean and of dimension $d_{j}<2n-1$. 
Then, the heat trace $\theta (\tau)$ defined by \eqref{theta} satisfies the condition {\rm (H)}, namely, 
for every sufficiently small $\ep>0$, 
there exists a constant $C>0$ such that, 
for any $\tau=x+iy$ with $x>0$ and $x \leq |y| \leq \pi$, we have 
\[
\re (\theta (x+iy)) -\theta (x) \leq -Cx^{-\ep}.  
\]
\end{prop}
\begin{proof}
Consider the distribution $\mu$ defined by 
\[
\mu=\sum_{\ell=1}^{\infty}e^{-i\lambda_{\ell}t}=\tr (e^{-itP})=\theta(it).
\]
Then, it is well-known (\cite{Ch}, \cite{DG}) that the singular support of $\mu$ is 
contained in the set of periods of the Hamilton flow of $p$. 
Let $T \in [-\pi,\pi]$. If $T$ is not a period, then $\mu$ is smooth 
around $T$. 
Since the distribution $\mu$ is the boundary value of the function 
$\theta (x+iy)$ as $x \searrow 0$, and since the function $\theta (\tau)$ is holomorphic 
on $\re (\tau)>0$, $\theta (x+iy)$ is bounded as $x \searrow 0$ when $y$ is close to
the regular point $T$. 
Now, suppose that $0 \neq T \in [-\pi,\pi]$ is a period of the Hamilton flow. 
Then, according to Theorem 4.5 in \cite{DG}, there exists an interval around $T$ 
in which no other periods occur, and on such an interval, there exists functions $\alpha_{1}(s),\ldots,\alpha_{r}(s)$ such that 
\[
\mu=\sum_{j=1}^{r}\int_{-\infty}^{+\infty}e^{-ist}e^{iTs}\alpha_{j}(s)\,ds 
\]
(the Fourier transform of $e^{iT\cdot}\alpha_{j}$ in the distribution sense) 
and the functions $\alpha_{j}(s)$ admits an asymptotic expansion of the form: 
\[
\alpha_{j}(s)\sim 
\sum_{k=0}^{\infty}\alpha_{j,k}s^{(d_{j}-1)/2-k},\quad s \to +\infty, 
\]
where $d_{j}$ is the dimension of the submanifold $Z_{j}$ and it is assumed that $d_{j} \leq 2n-2$. 
The functions $\alpha_{j}(s)$ are of order $|s|^{-N}$ as $s \to -\infty$ for any $N>0$.

Then, we note that the function $\theta(\tau)$ with $\re (\tau)>0$ can be written as 
\[
\theta (\tau) =\sum_{j=1}^{r}\int e^{-\tau s}e^{iTs}\alpha_{j}(s)\,ds. 
\]
In fact, the right hand side of the above, which we denote $H(\tau)$, 
is holomorphic for $\re (\tau) >0$ whose boundary value as a distribution is $\mu$, 
and which equals that of $\theta (\tau)$. 
Thus, since the boundary value of the holomorphic function $\theta -H$ is zero, 
we have $\theta=H$ near $T$ (see, for example, \cite{Ho2}). 
We then introduce the cut-off function $\chi$ such that $\chi =1$ on $[1,\infty)$ 
and $\chi =0$ on $(-\infty,0]$, and we set 
\[
\begin{gathered}
\theta (\tau)=\sum_{j=1}^{r}[F_{j}(\tau) +G_{j}(\tau)],\\
F_{j}(\tau)=\int e^{-\tau s}\chi(s)e^{iTs}\alpha_{j}(s)\,ds,\quad 
G_{j}(\tau)=\int e^{-\tau s}(1-\chi(s))e^{iTs}\alpha_{j}(s)\,ds.  
\end{gathered}
\]
Clearly, $|G_{j}(\tau)|$ is bounded when $x \searrow 0$ and $y$ is close to $T$, 
where we denote $\tau =x+iy$. 
Changing the variable $t=xs$ and using the estimate $|\alpha_{j}(s)| \leq cs^{(d_{j}-1)/2}$, we have 
\[
|F_{j}(\tau)| \leq cx^{-(d_{j}-1)/2-1}. 
\]
Since we assumed that $d_{j} \leq 2n-2$, we have $x^{-(d_{j}-1)/2-1} \leq cx^{1/2-n}$. 
Therefore, we have 
\[
\re (\theta (x+iy)) -\theta (x) \leq 
|\theta (x+iy)|-\theta (x) \leq -c_{0}x^{-n}+cx^{1/2-n} +C
\]
for some constants $c_{0},c,C>0$. Thus, the estimate in the assertion 
follows when $y$ is close to the period $T \neq 0$.

It is rather easy to handle the case where $y$ is close to zero. 
Note that, by the ellipticity of the operator $P$, 
the function $\Pi(x,\xi)$ on $(x,\xi) \in T^{*}M \setminus 0$ defined by 
\[
\Pi(x,\xi)=\inf\{t>0\,;\,\varphi_{t}(x,\xi)=(x,\xi)\}
\]
is lower semi-continuous strictly positive function on $\Sigma$. 
Since $\Sigma$ is compact, the function $\Pi$ has 
a minimum $t_{0} > 0$ on $\Sigma$. Thus, there is an interval around zero 
which does not contain any positive period. 
Thus, we can express $\theta (\tau)$ as 
\[
\theta (\tau)=c_{0}\tau^{-n} +a(\tau),\quad 
|a(\tau)| =O(|\tau|^{1-n}) 
\]
with a positive constant $c_{0}$ (see the proof of Corollary $2.2'$ in \cite{DG}). 
Then, by \eqref{aux11}, we have 
\[
\re (\theta(x+iy)) -\theta (x) \leq 
-c_{1}x^{-n}(1+o(x)), 
\]
which shows the asserted estimate. 
\end{proof}

\begin{prop}
\label{errorE1}
Fix a real number $\beta$ such that $1 < \beta <1+\frac{n}{2}$. 
Assume that the heat trace $\theta (\tau)$ defined by \eqref{theta} satisfies the condition {\rm (H)} 
$($which is explicitly stated in Proposition $\ref{condH}$$)$. 
Then, there exists constants $c>0$, $C>0$ and $\ep_{n}>0$ such that 
\[
\left|
G_{P}(x+iy)\exp [-\sum_{j=0}^{n-1}K_{j}x^{-(n-j)}]
\right| 
\leq Ce^{-cx^{-\ep_{n}}}
\]
for $x+iy$ with $x>0$, $x^{\beta} \leq |y| \leq \pi$. 
\end{prop}
\begin{proof}
As in \cite{Me}, \cite{An}, we devide the region $x^{\beta} \leq |y| \leq \pi$ into two parts. 
Note that, for the proof of the main theorem, 
we need the case that $0<x$ is sufficiently small. So, we may devide it into 
$x^{\beta} \leq |y| \leq x$ and $x \leq |y| \leq \pi$. 
First, we consider the region $x^{\beta} \leq |y| \leq x$. 
By a direct computation using Lemma \ref{formula1}, we have 
\[
\begin{split}
A(x+iy)&:=
\left|
G_{P}(x+iy)\exp [-\sum_{j=0}^{n-1}K_{j}x^{-(n-j)}]
\right| \\
&\leq C\exp [
\sum_{j=0}^{n-1}K_{j}
\{\re (x+iy)^{-(n-j)}-x^{-(n-j)}\}
+O(|\log x|)
],
\end{split}
\]
where $C>0$ is a constant. 
We set $\omega =y/x$. Then, $|\omega| \leq 1$. 
It is easy to show that 
there exists constants $c_{1}, c_{2}>0$ such that 
\[
-c_{1} x^{-(n-j)} \omega^{2} \leq \re (x+iy)^{-(n-j)} -x^{-(n-j)} \leq -c_{2}x^{-(n-j)}\omega^{2}. 
\]
Therefore, we get 
\[
A(x+iy)  \leq C\exp
\left[
-c_{3}x^{-n}\omega^{2}+c_{4}\omega^{2}\sum_{j=1}^{n-1}x^{-(n-j)} +O(|\log x|)
\right]. 
\]
Now, we have $x^{2(\beta -1)} \leq \omega^{2} \leq 1$ and $0<n-2(\beta -1)$. Thus we have 
\[
A(x+iy) \leq C\exp[-c_{5}\omega^{2}x^{-n}] \leq C\exp [-c_{5}x^{-n+2(\beta -1)}]. 
\]
Therefore, setting $\ep_{n}=n-2(\beta -1)$, we obtain the desired inequality. 
Next, consider the case where $x \leq |y| \leq \pi$. 
First, we note that $\re \log (G_{P}(x+iy)) =\log |G_{P}(x+iy)|$. 
By using \eqref{logT} for $\tau=x+iy$, we have 
\[
-\re \log(1-e^{-\lambda_{\ell}(x+iy)}) \leq e^{-\lambda_{\ell}x}\cos (\lambda_{\ell}y)-\log (1-e^{-\lambda_{\ell}x})-e^{-\lambda_{\ell}x}.  
\]
Summing over $\ell \geq 1$ in the above, we have 
\[
\log |G_{P}(x+iy)| 
\leq \re \theta (x+iy)+\log G_{P}(x) -\theta (x). 
\]
Since $x \leq |y| \leq \pi$, we can use the condition (H), that is, for each sufficiently small $\varepsilon>0$, 
there exists a positive constant $C$ such that, for any $\tau =x+iy$ with $x>0$, $x \leq |y| \leq \pi$, we have 
\[
\re (\theta (x+iy)) -\theta (x) \leq -Cx^{-\varepsilon}. 
\]
Hence we have 
\[
\log |G_{P}(x+iy)|
\leq \log G_{P}(x) -Cx^{-\varepsilon}. 
\]
Therefore, applying Lemma \ref{formula1} for $\tau=x$ and taking $\varepsilon < n-2(\beta -1)$ small enough, 
we have 
\[
A(x+iy) \leq \exp[-Cx^{-\varepsilon} +O(|\log x|)], 
\]
which shows the assertion. 
\end{proof}

\section{Proof of the main theorem}
\label{Proof}

The purpose of this section is to give proofs 
of the main theorems. First of all, we will give a definition 
of the positive number $x_{E}$ appeared in the statement 
of Theorem \ref{main1}. 

\subsection{The positive number $x_{E}$}
\label{Proof1}

Define the polynomial $p_{E}(x)$ with one variable by 
\begin{equation}
\label{spectralP}
p_{E}(x)=Ex^{n+1}-\sum_{j=0}^{n-1}
(n-j)K_{j}x^{j}, 
\end{equation}
where the constants $K_{j}$ are defined in \eqref{residue} by 
using the residues of poles of the spectral zeta function. 
The polynomial $p_{E}$ is defined by the spectral data and 
the fixed energy level $E$. 
Since $K_{0}=A_{0}\zeta(n+1)\Gamma (n)$, $A_{0}=(2\pi)^{-n}\vol(\Sigma)$ is positive, 
we have $p_{E}(0)=-nK_{0} <0$. The coefficient of $x^{n+1}$ is $E>0$, and hence 
we have $p_{E}(x) \to +\infty$ as $x \to \infty$. Therefore, 
there is a solution of the equation $p_{E}(x)=0$ in the positive real axis. 
Now, we have the following 
\begin{prop}
For sufficiently large $E>0$, the polynomial $p_{E}(x)$ has just one positive zero. 
We denote the zero of $p_{E}$ by $x_{E}$ {\rm (}as in Theorem {\rm \ref{main1}}{\rm )}. 
Then, the asymptotics of $x_{E}$ as $E \to \infty$ is given by \eqref{zeroE}.
\end{prop}
\begin{proof}
First of all, let $y_{E}>0$ be a positive zero of the polynomial $p_{E}$. 
Then, we have 
\[
0< E  \leq \sum_{j=0}^{n-1}(n-j)|K_{j}|y_{E}^{j-n-1}. 
\]
It follows from this that $y_{E}=O(E^{-\frac{1}{n+1}})$ as $E \to \infty$. 
We have 
\[
p'_{E}(y_{E})=n(n+1)K_{0}y_{E}^{-1}
+\sum_{j=1}^{n-1}(n-j)(n+1-j)K_{j}y_{E}^{j-1}. 
\]
The sum in the right hand side of the above is bounded as $E \to \infty$, 
and the first term goes to $\infty$. (Each $K_{j}$ is a real number.) 
Therefore, we have $p'_{E}(y_{E}) >0$ 
for each zero $y_{E}$ of $p_{E}$ for sufficiently large $E>0$. 
Thus, for sufficiently large $E>0$, the polynomial $p_{E}$ has 
a unique positive zero $x_{E}>0$. 
We note that 
\begin{equation}
\label{asyX}
Ex_{E}^{n+1}=nK_{0}+\sum_{j=1}^{n-1}(n-j)x_{E}^{j}=nK_{0}+O(E^{-\frac{1}{n+1}}), 
\end{equation}
which shows the asymptotics \eqref{zeroE}.
\end{proof}

We notice that, without any assumption on the the eigenvalues of $P$ and classical Hamilton flow, 
we can prove the upper bound for $\Omega$ by using Lemma \ref{formula1} 
which agree with the limit formula in Corollary \ref{limL}. 

\begin{prop}
\label{upperBB}
We have 
\begin{equation}
\label{upperB}
\Omega (E) \leq c E^{\frac{Z_{P}(0)}{n+1}}e^{B_{n}E^{\frac{n}{n+1}}(1+o(1))},\quad 
B_{n}=(n+1)\left(
\frac{\vol (\Sigma)}{(2\pi n)^{n}}
\zeta (n+1)\Gamma (n)
\right)^{\frac{1}{n+1}}
\end{equation}
as $E \to \infty$ with a positive constant $c$. 
\end{prop}
\begin{proof}
We set 
\[
\psi (x)=\sum_{j=0}^{n-1}K_{j}x^{-(n-j)},\quad x >0. 
\]
We take arbitrary $M>0$. Then, by Lemma \ref{formula1}, we can find positive 
constants $c_{1}$, $c_{2}$ (depending on $M$) such that 
\[
c_{1}x^{-a}e^{\psi(x)} \leq G_{P}(x) \leq c_{2}x^{-a}e^{\psi (x)}
\]
for any $x$ in the interval $(0,M]$, where we set $a=Z_{P}(0)$. Thus, for any 
$E \in \Gamma (P)$ and $0<x \leq M$, we have 
\[
\Omega(E) \leq c_{2}x^{-a}e^{\psi (x)+Ex}. 
\]
We take a sufficiently large positive number $E_{o}$ such that $x_{E} <M$ for any $E \geq E_{o}$, 
and, for $E \geq E_{o}$, we set $x =x_{E}$, so that the function $\psi (x) +Ex$ takes 
its minimum there. Then, a simple computation with \eqref{zeroE} shows 
\[
\psi (x_{E}) +Ex_{E}=B_{n}E^{\frac{n}{n+1}}(1+O(E^{-\frac{1}{n+1}})), 
\]
and which completes the proof. 
\end{proof}

\subsection{Proof of Theorem \ref{WeylAA}}
We use Hardy-Ramanujan Tauberian theorem (\cite{HR1}, Theorem A) 
which asserts the following: Suppose the sequences $\{E_{m}\}_{m=1}^{\infty}$, $\{a_{m}\}_{m=1}^{\infty}$ of 
non-negative real numbers and the constants $A>0$, $\alpha>0$ satisfy the following: 
\begin{enumerate}
\item $E_{m}>E_{m-1}$, $E_{m} \to \infty$;
\item $E_{m}/E_{m-1} \to 1$;
\item $G(s):=\sum a_{m}e^{-E_{m}s}$ is convergent for $s>0$; 
\item $G(s)=\exp [As^{-\alpha}(\log (1/s))^{-\beta}(1+o(1))]$ as $s \to 0$, where $\beta$ is a real number. 
\end{enumerate}
Then we have 
\[
A_{n}:=\sum_{m=1}^{n}a_{m}=\exp[BE_{n}^{\alpha/(1+\alpha)}(\log E_{n})^{-\beta/(1+\alpha)}(1+o(1))],
\]
where $B$ is given by 
\[
B=A^{1/(1+\alpha)}\alpha^{-\alpha/(1+\alpha)}(1+\alpha)^{1+\beta/(1+\alpha)}. 
\]
In our situation, we denote $\Gamma (P)=\{0=E_{0} <E_{1}<E_{2}<\cdots\}$.  
Then we have $E_{m}-E_{m-1} \leq \lambda_{1}$, 
and hence $E_{m}/E_{m-1} \to 1$ as $m \to \infty$. 
Thus the condition $(1)$ and $(2)$ are satisfied. The non-negative number $a_{m}$ 
in the condition $(3)$ corresponds to the number of states $\Omega(E_{m})$ of Boson gas, 
and hence $(3)$ follows from Lemma \ref{pfunc}. The condition $(4)$ is deduced from 
Lemma \ref{formula1} with $\alpha=n$, $\beta=0$ and $A=K_{0}=(2\pi)^{-n}\vol (\Sigma)\zeta(n+1)\Gamma(n)$. 
Thus, the formula \eqref{WeylA} follows from the Hardy-Ramanujan Tauberian theorem with $B=B_{n}$.

\subsection{Proof of Theorem \ref{main1}}
\label{Proof2}

First of all, we fix $0<\mu <\frac{1}{n}$ ($0<\mu<1/2$ for $n=1$) 
and choose $0<C_{0}<1$, $1<\beta <1+\frac{n}{2}$ as follows. 
For $n=1$, we set $\beta=\frac{11}{8}+\frac{\mu}{4}$ and $C_{0}=\frac{1}{8}+\frac{3}{4}\mu$. 
For $n = 2$, we set $\beta=\frac{7}{4}+\frac{1}{2}\mu$ and $C_{0}=\frac{1}{4}+\frac{3}{2}\mu$. 
For $n=1,2$, we set $\delta=\frac{1}{2}-\mu$. 
For $n \geq 3$, we take $\delta$ satisfying $0<\delta < \min\{\frac{1}{2}, \frac{4}{n}, \frac{4}{3}(\frac{1}{2}-\frac{1}{n}), 4(\frac{1}{n}-\mu)\}$ 
and set $\beta=1+\frac{n}{2}(1-\frac{\delta}{2})$ and $C_{0}=n(\mu+\frac{\delta}{4})$. 
In all cases, the positive numbers $\beta$, $\delta$ and $C_{0}$ satisfy $0<\delta<\min\{1/2,4/n\}$, 
$1+\frac{n}{3} < \beta=1+\frac{n}{2}(1-\frac{\delta}{2}) <1 +\frac{n}{2}$, $0<C_{0}<1$ and $\mu=\min\{\frac{1}{2}-\delta, \frac{C_{0}}{n}-\frac{\delta}{4}\}$.

From now on, we assume that the eigenvalues of $P$ are integers. 
Then, we have $\Gamma (P) \subset \mb{Z}$, and hence we can use the integral formula \eqref{int2} 
for the number of states $\Omega(E)$. 
For every $x>0$, we have 
\begin{equation}
\label{se1}
\begin{gathered}
\Omega(E)=I_{1}(E)+R_{1}(E), \\
I_{1}(E)=\frac{1}{2\pi}\int_{|y| \leq x^{\beta}}e^{E(x+iy)}G_{P}(x+iy)\,dy, \\ 
R_{1}(E)=\frac{1}{2\pi}\int_{x^{\beta} \leq |y| \leq \pi}e^{E(x+iy)}G_{P}(x+iy)\,dy, 
\end{gathered}
\end{equation}
where $\beta$ is a fixed constant as above. 
We shall determine $x>0$ by using the suddle point argument 
for the integral $I_{1}(E)$. By Lemma \ref{formula1}, we have 
\[
I_{1}(E)=\frac{1}{2\pi}\int_{|y| \leq x^{\beta}}
e^{E\tau+\sum_{j=0}^{n-1}K_{j}\tau^{j-n}}e^{Z_{P}'(0)-Z_{P}(0)\log \tau +J(\tau)}\,dy,
\quad \tau=x+iy. 
\]
Then, in the integrand of the above, we think $E\tau +\sum_{j=0}^{n-1}K_{j}\tau^{j-n}$ 
as a phase function. The critical point of this phase function close to $\im (\tau)=0$ 
is given by $\tau =x_{E}$. Thus, we choose $x=x_{E}$ in \eqref{se1}. 
First, let us consider the integral $I_{1}(E)$ in \eqref{se1}. 
As in \cite{An}, we write $m_{E}=Ex_{E}$. 
Then, by Lemma \ref{formula1} and the change of variable $\omega=y/x_{E}$, we have 
\[
\begin{split}
I_{1}(E)&=
\frac{1}{2\pi}e^{
m_{E}+Z_{P}'(0)+(1-Z_{P}(0))\log x_{E}
}\\
& \times \int_{|\omega|\leq x_{E}^{\beta-1}}
e^{
im_{E}\omega +
\sum_{j=0}^{n-1}K_{j}x_{E}^{-(n-j)}(1+i\omega)^{-(n-j)}
-Z_{P}(0)\log (1+i\omega)
+O(m_{E}^{-C_{0}/n})
}\,d\omega,
\end{split}
\]
where we have used $x_{E}=O(E^{-\frac{1}{n+1}})=O(m_{E}^{-1/n})$. 
By using the formula 
$(1+i\omega)^{-\alpha}=1-i\alpha \omega 
-\frac{1}{2}\alpha(\alpha +1)\omega^{2}
+O(|\omega|^{3})$, and 
\[
m_{E}=Ex_{E}=\sum_{j=0}^{n-1}(n-j)K_{j}x_{E}^{-(n-j)}, 
\]
we obtain 
\begin{equation}
\label{se2}
\begin{gathered}
I_{1}(E)=\frac{1}{2\pi}e^{
x_{E}^{-n}f(x_{E})+Z_{P}'(0)+(1-Z_{P}(0))\log x_{E}
}\times I_{2}(E), \\
I_{2}(E)=
\int_{|\omega| \leq x_{E}^{\beta -1}}
e^{
-\frac{1}{2}\eta_{n}(E)\omega^{2}
-Z_{P}(0)\log (1+i\omega) +O(m_{E}^{-C_{0}/n})+
\sum_{j=0}^{n-1}K_{j}x_{E}^{-(n-j)}O(|\omega|^{3})
}\,d\omega, \\
\eta_{n}(E)=
\sum_{j=0}^{n-1}
(n-j)(n-j+1)K_{j}x_{E}^{-(n-j)}, 
\end{gathered}
\end{equation}
where the function $f$ is defined in \eqref{freqF}. By the definition of 
the number $x_{E}$ and the function $f$, one has 
\[
x_{E}^{-n}f(x_{E})=Ex_{E}+\sum_{j=0}^{n-1}K_{j}x_{E}^{-(n-j)}. 
\]
Now, it is easy to see that 
$x_{E}^{-(n-j)}|\omega|^{3}=O(m_{E}^{1-\frac{3}{n}(\beta -1)})$ 
for every $j=0,\ldots,n-1$. 
We note that $1-\frac{3}{n}(\beta -1) <0$. 
Near $\omega=0$, obviously we have $\log (1+i\omega)=O(|\omega|)$, 
and, since $\beta <1+\frac{n}{2}$, we have 
$|\omega|=O(m_{E}^{-\frac{1}{n}(\beta -1)}) =O(m_{E}^{1-\frac{3}{n}(\beta -1)})$. 
Thus, if we set 
$\mu_{1}=\min \{\frac{3}{n}(\beta -1)-1,\,C_{0}/n\}$, 
we obtain 
\[
I_{2}(E)=\int_{|\omega| \leq x_{E}^{\beta -1}}
e^{-\frac{1}{2}\eta_{n}(E)\omega^{2}+O(m_{E}^{-\mu_{1}})}\,d\omega
=A(E)+R_{2}(E), 
\]
where we set 
\[
\begin{gathered}
A(E)=\int_{|\omega| \leq x_{E}^{\beta -1}} 
e^{-\frac{1}{2}\eta_{n}(E)\omega^{2}}\,d\omega,\\
R_{2}(E)=\int_{|\omega| \leq x_{E}^{\beta -1}}
e^{-\frac{1}{2}\eta_{n}(E)\omega^{2}}
\left(
e^{O(m_{E}^{-\mu_{1}})}-1
\right)\,d\omega.
\end{gathered}
\]
By using the definition of $\eta_{n}(E)$ and \eqref{asyX}, we have 
\begin{equation}
\label{etaAsy}
m_{E}^{-1}\eta_{n}(E)=(n+1)+O(m_{E}^{-1/n}). 
\end{equation}
Thus, for sufficiently large $E>0$, $\eta_{n}(E) >0$. 
Therefore, we have 
\[
|R_{2}(E)|=O(m_{E}^{-\mu_{1}-(\beta -1)/n}). 
\]
Then, we easily have 
\[
A(E)=
\sqrt{\frac{2\pi}{\eta_{n}(E)}}+O(m_{E}^{-1/2}e^{-Cm_{E}^{1-2(\beta-1)/n}}). 
\]
For $\eta_{n}(E)$, by using \eqref{etaAsy}, we have 
\[
\frac{1}{\sqrt{\eta_{n}(E)}}=
\frac{1}{\sqrt{(n+1)m_{E}}}
(1+O(m_{E}^{-1/n})). 
\]
Therefore, since $\beta =1+\frac{n}{2}(1-\frac{\delta}{2})$ with 
$0 < \delta < \min \{1/2,4/n\}$ and 
\[
\mu=\min \{\frac{1}{2}-\delta,\,\frac{C_{0}}{n}-\frac{\delta}{4}\}, 
\]
we obtain 
\[
I_{2}(E)=\sqrt{\frac{2\pi}{(n+1)m_{E}}}(1+O(m_{E}^{-\mu})). 
\]
Inserting the definition $m_{E}=Ex_{E}$, and using \eqref{zeroE}, 
we obtain 
\[
I_{2}(E)=\sqrt{\frac{2\pi}{n+1}}
\left(\frac{1}{nK_{0}}\right)^{1/2(n+1)}E^{-\frac{n}{2(n+1)}}
(1+O(E^{-\frac{n\mu}{n+1}}))
\]
By \eqref{se2}, we conclude, for the integral $I_{1}(E)$, that 
\[
I_{1}(E)=\frac{e^{x_{E}^{-n}f(x_{E})}}{\sqrt{2\pi(n+1)}}
\left(
\frac{1}{nK_{0}}
\right)^{1/2(n+1)}
E^{-\frac{n}{2(n+1)}}e^{Z_{P}'(0)}x_{E}^{1-Z_{P}(0)}(1+O(E^{-\frac{n\mu}{n+1}})).
\]
Inserting, in the above, the formula 
\[
x_{E}^{1-Z_{P}(0)}=
\left(
\frac{nK_{0}}{E}
\right)^{\frac{1-Z_{P}(0)}{n+1}}
(1+O(E^{-\frac{1}{n+1}})),
\]
which is deduced from \eqref{zeroE}, we get the right hand side of \eqref{mainF}. 
Next, we have to consider the integral $R_{1}(E)$ in \eqref{se1}. 
By Proposition \ref{errorE1}, we have 
\[
|R_{1}(E)| \leq 
Ce^{-C_{n}x_{E}^{-\ep_{n}}}\exp
[Ex_{E}+\sum_{j=0}^{n-1}K_{j}x_{E}^{-(n-j)}]
=Ce^{-C_{n}x_{E}^{-\ep_{n}}}
e^{x_{E}^{-n}f(x_{E})}
\]
with positive constants $C$, $C_{n}$, $\ep_{n}$. 
Since the term $e^{-C_{n}x_{E}^{-\ep_{n}}}$ decays exponentially and 
since the term $e^{x_{E}^{-n}f(x_{E})}$ is the same as the exponential term in the asymptotics
of the integral $I_{1}(E)$, we conclude \eqref{mainF}. This completes the proof 
of Theorem \ref{main1}.

%

\end{document}